\documentclass[12pt]{article}
\usepackage{latexsym}
\begin{document}


\newtheorem{Thm}{\noindent Theorem}[section]
\newtheorem{Lem}{\noindent Lemma}[section]
\newtheorem{Cor}{\noindent Corollary}[section]
\newtheorem{Prop}{\noindent Proposition}[section]

\makeatletter
\renewcommand{\theequation}
{\thesection.\arabic{equation}}
\@addtoreset{equation}{section}




                

\setlength{\baselineskip}{15pt}
\newcommand{\vlimsup}{\mathop{\overline{\lim}}}
\newcommand{\vliminf}{\mathop{\underline{\lim}}}
\newcommand{\Av}{\mathop{\mbox{Av}}}
\newcommand{\spec}{{\rm spec}}
\newcommand{\qed}{\hfill $\Box$ \\}
\newcommand{\subarray}[2]{\stackrel{\scriptstyle #1}{#2}}
\def\textmc{\rm}
\def\({(\!(}
\def\){)\!)}
\def\R{{\bf R}}
\def\Z{{\bf Z}}
\def\N{{\bf N}}
\def\C{{\bf C}}
\def\T{{\bf T}}
\def\E{{\bf E}}
\def\H{{\bf H}}
\def\Prob{{\bf P}}
\def\M{{\cal M}}     
\def\F{{\cal F}}
\def\G{{\cal G}}
\def\D{{\cal D}}
\def\X{{\cal X}}
\def\A{{\cal A}}
\def\B{{\cal B}}
\def\L{{\cal L}}
\def\a{\alpha}
\def\b{\beta}
\def\e{\varepsilon}
\def\de{\delta}
\def\ga{\gamma}
\def\k{\kappa}
\def\la{\lambda}
\def\fa{\varphi}
\def\th{\theta}
\def\si{\sigma}
\def\t{\tau}
\def\om{\omega}
\def\De{\Delta}
\def\Ga{\Gamma}
\def\La{\Lambda}
\def\Om{\Omega}
\def\Th{\Theta}
\def\lan{\langle}
\def\ran{\rangle}
\def\lbr{\left(}
\def\rbr{\right)}
\def\const{\;\operatorname{const}} 
\def\dist{\operatorname{dist}} 
\def\Tr{\operatorname{Tr}}
\def\quadd{\qquad\qquad}
\def\n{\noindent}
\def\beq{\begin{eqnarray*}}
\def\eeq{\end{eqnarray*}}
\def\supp{\mbox{supp}}
\def\beqn{\begin{equation}}
\def\eeqn{\end{equation}}
\def\bp{{\bf p}}
\def\sg{{\rm sgn\,}}
\def\1{{\bf 1}}
\def\pf{{\it Proof.}}
\def\v2{\vskip2mm}
\def\n{\noindent}
\def\z{{\bf z}}
\def\x{{\bf x}}

\begin{flushleft}
{Available online at www sciencedirect. com }\\
{To appear in Stoch. Proc.  Appl.  {\bf 123} 191-211  (2013)}
\end{flushleft}
\begin{center}
\vskip10mm
{\Large  The expected area of  the Wiener sausage swept by a disc
 } \\
\vskip4mm
{K\^ohei UCHIYAMA} \\
\vskip2mm
{Department of Mathematics, Tokyo Institute of Technology} \\
{Oh-okayama, Meguro Tokyo 152-8551\\
e-mail: \,uchiyama@math.titech.ac.jp}
\end{center}

\vskip8mm
\n
{\it key words}: Wiener sausage, Brownian hitting time,  Brownian bridge, Ramanujan's function, asymptotic expansion,  Laplace inversion
\vskip2mm
\n

\vskip6mm

\begin{abstract}
The expected areas of the Wiener sausages  swept by a disc  attached to the two-dimensional Brownian Bridge  joining  the origin to a point $x$ over  a time interval $[0,t]$ are computed.  It is proved that  the leading term of the expectation is given by  Ramanujan's function if $|\x|=O(\sqrt t)$. The second term is also given explicitly  when $|\x|=o(\sqrt t)$. The corresponding result for unconditioned process is also obtained.
  \end{abstract}
\vskip6mm

\section { Introduction}
\vskip2mm

Let $B_t$ be the standard two-dimensional Brownian motion started at the origin defined on a probability space $(\Om,\F, P)$. Fixing $r>0$ let $S^{(r)}_t$ be the Wiener sausage of radius $r$ and length $t$, namely it is  the region  swept by the disc of radius $r$ attached to $B_s$ at its center as $s$ runs from $0$ to $t$:
$$S^{(r)}_t = \{\z\in \R^2: |B_s - \z| <r~~\mbox{for some} ~~s\in [0,t] \}.$$
In this paper we compute the expectation of the  area of $S^{(r)}_t $, which we denote by 
${\rm Area} (S^{(r)}_t)$, for  Brownian motion  conditioned  to be at a prescribed point  $\x \in \R^2$ at time $t$ as well as for  free (unconditioned) Brownian motion.  Define $N(\la)$, called Ramanujan's function (\cite{B}, \cite{W}) or integral (\cite{E2}, page 219), by   
$$N(\la) = \int_0^\infty \frac{e^{-\la u}}{(\lg u)^2 +\pi^2}\cdot\frac{du}{u}~~~~~~~(\la \geq 0).$$
Put $\k = 2e^{-2\ga }$ where $\ga = -\int_0^\infty e^{-u}\lg u\,du$ (Euler's constant).  
 Let $E$ designate  the expectation with respect to $P$.   In this paper we prove the following two theorems.

\vskip2mm 
\begin{Thm}\label{thm1} ~ 
\beq
\frac{d}{dt} E[{\rm Area} (S^{(r)}_t)]&=& 2\pi   N\bigg(\frac{ \k t}{r^2}\bigg)  + \frac{4\pi r^2}{t(\lg (t/r^2))^3}(1+o(1))   ~~~~~~~~~\mbox{as}~~~~~t   \to   \infty,\\
&=&r  \sqrt {\frac{2\pi }{ t}} + \frac{\pi}{2}  + O\Big( \sqrt t \,\Big)  ~~~~~~~~~~~~~~~~~~~~~~~~~~\mbox{as}~~~~~t \to 0.
\eeq
\end{Thm}

For $\x\in \R^2$ we write $\x^2$ for the square of Euclidian length $|\x|$.
  \vskip2mm 
  \begin{Thm}\label{thm2} ~ For each $M>1$, uniformly for $|\x|<M \sqrt{t}$, as $t\to \infty$ 
  $$E[{\rm Area} (S^{(r)}_t) \,|\, B_t=\x] = 2\pi t N\bigg(\frac{\k t}{r^2}\bigg) + \frac{\pi \x^2}{(\lg t)^2}\bigg[\lg\bigg(\frac{t}{\x^2\vee 1}\bigg) +O(1)\bigg]  +O(1), $$
  where the left-hand side is the conditional expectation conditioned on $B_t=\x$.
\end{Thm}

      The proof of Theorem \ref{thm1} is performed  by Laplace inversion with rather simple computations.  The leading term in the formula of Theorem 
      \ref{thm2}  is derived in a similar way to  that in  Theorem \ref{thm1}, while  the identification of the second term of it  is made by  a delicate analysis unless $|\x|$ remains in a bounded set. 
  \vskip2mm\n

{\sc Remark 1.}~ 
The function $N(t)$  admits the following asymptotic expansion in powers of $1/\lg t$ as $t\to\infty$:
\begin{eqnarray} \label{Rmk1}
N( t) \sim \frac1{\lg t}+\frac{-\ga }{(\lg t)^{2}}+\frac{\,\ga ^2- \zeta(2)\,}{(\lg t)^{3}}
+    \frac{\,3\ga \zeta(2)+ (-\ga)^3 - 2\zeta(3)}{(\lg t)^{4}\,} + \cdots,
\end{eqnarray}
 where $\zeta(z)= \sum_{n=1}^\infty n^{-z}$. While it requires a tedious computation to derive  the corresponding expansion of $N(\k t/r^2)$ directly from this one, actually there is a simple way to transform the expansion.  For each $\a>0$ the expansion for $N(\a t)$ is obtained by simply replacing $-\ga$ by $-\ga-\lg \a$ in (\ref{Rmk1}), so that the expansion up to the third order term becomes
\begin{eqnarray} \label{Rmk2}
N(\a t) \sim \frac1{\lg t}+\frac{-\ga -\lg \a}{(\lg t)^{2}}+ \frac{(\ga +\lg \a)^2- \zeta(2)}{(\lg t)^{3}}
+ \cdots 
\end{eqnarray}
Similarly the asymptotic expansion of $\frac1{t}\int_0^t N(\a s)ds$ is obtained by replacing $-\ga$ by $1-\ga-\lg \a$ and $-\zeta(k)$  by $1-\zeta(k)$ in (\ref{Rmk1}),
so that
 \begin{eqnarray} \label{Rmk3}
\frac1{t}\int_0^t N(\a s)ds \sim
\frac{1}{\lg t}+\frac{1-\ga -\lg \a}{(\lg t)^{2}}+\frac{(1-\ga -\lg \a)^2  +1- \zeta(2)}{(\lg t)^{3}}  
+ \cdots  
\end{eqnarray}
(See Section 5 for  these and related  matters.)    It also is noted that  combining  Theorems \ref{thm1} and \ref{thm2}  yields
  \beq
    E[{\rm Area} (S^{(r)}_t) ] -E[{\rm Area} (S^{(r)}_t) \,|\, B_t=0] 
  = 2 \pi    \int_0^{ t}\Big[ -\a sN'(\a s)\Big]ds 
 +O(1)
  \eeq
as $t\to\infty$, where $\a =\k /r^2$.  The asymptotic expansion of this difference is obtained from   that of $ E[{\rm Area} (S^{(r)}_t) ]/\lg t$  by a very simple rule (see Lemma \ref{lem_RF2} in Section 5).

\v2
{\sc Remark 2.}~ It may be reasonable to compare the expected increment of the sausage at time $t$ with $2\pi r E[|B_t|]$ as $t\downarrow 0$. 
The former one is asymptotic to $2r \sqrt{2\pi t}$ according to  the second half of Theorem \ref{thm1}, while the latter equals  $2\pi r \sqrt{\pi t/2}$, so that
the ratio of the latter to the former equals $\pi/2 $.
\v2

  Asymptotic behavior of the area ${\rm Area} (S^{(r)}_t)$  (or volume in higher dimensions) as $t\to\infty$   has long been studied from various points of view. It
    is a typical functional  of Brownian paths that is non-Markovian and the standard limit theorems for it  (the law of large numbers, the central limit theorems or the large deviations) have been  of continued interest (\cite{Wh},  \cite{LG2}, \cite{HK},  etc.) .
  The expectation of it is the total heat emitted in the  time interval  $[0,t]$  from the disc which is kept  at the unit temperature. The sausage for conditioned process (with $\x=0$)  naturally arises in the study of the asymptotic estimate of the trace of the heat kernel on the plane with randomly scattered cooling discs kept at zero temperature and has been effectively used (cf. \cite{KL}, \cite{DV}). 
 
 For  free Brownian motion, some asymptotic  expansions in negative powers of $\lg t$ of the expectation divided by $t$  for the sausage swept by an arbitrary (non-polar) compact  subset of $\R^2$  are obtained by Spitzer \cite{S2} up to magnitude $o(1/(\lg t)^2)$  and by   Le Gall \cite{LG} to any order (see also \cite{LG3}).   
M. van den Berg and E. Bolthausen  \cite{BB} compute the expectation for  Brownian motion conditioned on returning to the origin at time $t$ for the disc case  up to the error term of magnitude
 $O(t\sqrt{\lg \lg t)}/ (\lg t)^4$ and conjecture that  their formula for the conditional  expansion would be valid for the general compact set, $K$ say, of  positive capacity if the radius of the disc is  replaced by the logarithmic capacity of $K$ (in the usual normalization \cite{Af}), the situation already observed for the free Brownian case.   
 
 The  conjecture stated in \cite{BB} has been verified   in a very recent paper \cite{M2}  by I. McGillivray:
 in fact he obtains  the asymptotic expansion for any compact set  $K$ of positive capacity and computes the explicit forms of the first three coefficients of the  expansion in terms of  the logarithmic capacity of $K$, the coefficients agreeing with those conjectured  in \cite{BB}.  
 It is warned  that these authors define the sausage for the Brownian motion $\tilde B_t = B_{2t}$ instead of $B_t$ so that
to translate our results  to their case  one must replace $t$  by $2t$ in our formulae. 
The coefficients of our formulae of course agree with those obtained previously,  whether it is a  free  or conditioned Brownian motion  which  forms the  sausage,  as readily ascertained  by substituting  $\a =\k/r^2$ in 
(\ref{Rmk2}) and (\ref{Rmk3}) and noting that the logarithmic capacity of the disc of radius $r$ equals $r^2$ (according to the normalization of logarithmic capacity in them).  This  suggests that our formulae in Theorems \ref{thm1} and \ref{thm2} would be extended to the sausage swept by any non-polar compact set  $K$ in place of the  disc if $r$ is replaced by the logarithmic capacity of $K$.

The higher dimensional case of  free Brownian motion is treated by  Spitzer \cite{S2}, \cite{G} and Le Gall \cite{LG}  for a compact set and by \cite{H} for a ball. The pinned cases are dealt with by Uhlenbeck and  Beth \cite{Uhl} and McGillivray \cite{M}. A more detailed account of the results for the expected volume (for dimensions $\geq 2$)  obtained  up to 1997 can be found in \cite{BB}.

The function $N(t)$  expressing the  leading terms in our theorems has already  appeared  in an analogous manner for the corresponding problem concerning  the range of random walks (cf. \cite{U1}) (in which $N'(\la)$ is denoted by $W(\la)$).  
The expression by means of $N(t)$ makes possible a subtle  comparison between the expected area of the Wiener sausage and the corresponding quantity to  a random walk of mean zero:  we can  associate a certain natural \lq radius', say $r_*$, of a lattice point of  $\Z^2$  with the  random walk (\cite{U2}, Remark 6),  and  according to Corollary 1.1 of  \cite{U1} and our Theorem \ref{thm2}  the expected number of sites  visited by the walk in the first  $n$ steps and  the expected  area of the Wiener  sausage over the internal  $[0,n]$ coincide up to 
the error of $O(1)$ for the processes  conditioned to return to  the origin at the  time $n$,  provided that $r$ is chosen to be $r_*$,  the variances of Brownian and random walk processes are the same and the fourth moments of the random walks exist.

The results obtained in this paper are quite parallel with those in \cite{U1}, but the methods and  the structures of the proofs   are considerably  different.  In the random walk  case
we have simple expressions of the expected \lq area' of the range of the walk for both free and conditional ones due to the discrete nature of the walk, which is not available to the Brownian case. In (cf. \cite{U1}) the Fourier analysis is the main tool, while in the present paper it plays a minor,  though fundamental, role.  We need a careful evaluation of certain integrals    to find out  the form of the second term  as given in the formula of Theorem \ref{thm2}. Our derivation of it  is directed by the result of \cite{U1}, the agreement of  asymptotic forms of various quantities for Brownian motion and Random walks being  expected. 
The  precise estimates of the  first hitting time distributions for corresponding processes as obtained in \cite{U2} and \cite{Ubh}   are used in both papers but the usage is much more essential for the present  than in \cite{U2}.  
One can obtain  an asymptotic estimate of the density of the  first hitting time distribution valid uniformly with respect to  starting positions and, with it,   extend the result of Theorem \ref{thm2} to the case when $|\x|/\sqrt t \to \infty$,  which  will be studied in  another paper \cite{Usaus2}.
   
   The following notation will  be used:  $a\wedge b =\min\{a,b\}$, $a\vee b=\max \{a,b\}$ ($a,b\in \R)$; two dimensional points are denoted by bold face letters $\z, \x, {\bf y}$,  $|\z|$ denotes the Euclidean length of $\z$; and $\x\cdot \z$ the Euclidian inner product of $\x$ and $\z$; for functions $g$ and $G$ of a variable $x$,
 $g(x)=O(G(x))$  means  that there exists  a constant $C$ such that $|g(x)|\leq C|G(x)|$ whenever  $x$ ranges over a specified set; the letters
 $C, C', C''$ etc. denotes constants whose values are not significant and  may change in different places where they occur.
 
 We prove  Theorem \ref{thm1}   in Section 2 and   Theorem \ref{thm2}  in Section 4. In Section 3 we give some results on the distribution of the first hitting time to a disc, which prepare for the analysis made in Section 4. In the last section we derive  the asymptotic expansions associated with $N(t)$ as mentioned in Remark 1.

\section{Proof of Theorem \ref{thm1}}  

We shall consider the Brownian motion started at $\z$ and denotes its law by $P_\z$.
 Let $\sigma^{(r)} =\sigma_{U(r)}$ be the first hitting  time of  Brownian motion $B_t$ to the  disc $U(r)$ of radius $r$  and centered at the origin. Then
  \beq
E[{\rm Area} (S^{(r)}_t)]  = \int_{\R^2} P_\z[\sigma^{(r)}<t\,] d \z.
\eeq
Let $q(\z;r)$ denote the density of the distribution of $\sigma^{(r)}$: 
$$q(\z, t;r)=\frac{d}{dt} P_\z[\sigma^{(r)}<t\,],$$ 
so that
\beqn\label{2.21}
\frac{d}{dt} E[{\rm Area} (S^{(r)}_t)]  =\int_{\z>r} q(\z, t;r) d\z.
\eeqn
Let  $p_t(\z)$ denote the heat kernel on the plane.:
$$p_t(\z)=(2\pi t)^{-1}e^{- \z^2/2t}.$$
As is well-known we have
$$\int_0^\infty p_t(\z)e^{-\la t}dt =\frac1{\pi} K_0(|\z|\sqrt{2\la})$$
(cf. \cite{IM}, \cite{E}), hence
\beqn\label{2.22}
\int_0^\infty q(\z,t;r)e^{-\la t}dt = \frac{G_\la(0,|\z|)}{G_\la(0,  r)} = \frac{ K_0(|\z|\sqrt{2\la})}{ K_0(r\sqrt{2\la})}, 
\eeqn
where   $G_\la$ denotes the resolvent kernel for the 2-dimensional Bessel process and $ K_{\nu}$  the usual modified Bessel function of second kind of order  $\nu$. 
Put
$$m(t;r) =\int_{|\z|>r} q(\z, t;r)d\z.$$
From (\ref{2.21}) and (\ref{2.22}) we have
 \begin{eqnarray}\label{2.23}
 \int_0^\infty m(t;r)e^{-\la t}dt 
 &=& 2\pi \int^\infty_{r} \frac{K_0(u\sqrt{2\la})}{K_0(r \sqrt{2\la})}\,u du \nonumber\\
 &=&  \frac{2\pi r K_1(r \sqrt{2\la})}{K_0(r \sqrt{2\la}) \sqrt{2\la}\,},
 \end{eqnarray}
 where we have applied  the identity $(d/dz)[z K_1(z)] = - zK_0(z)$ for the second equality. 
 From the scaling property of Brownian motion it follows that
 $$m(t,r) =m(t/r^2;1)$$
and   we suppose $r =1$ in what follows. By   Laplace inversion 
\beqn\label{BE0}
m(t, 1) = \int_{-\infty}^{\infty}   \frac{K_1( \sqrt{2iu})}{K_0( \sqrt{2iu}) \sqrt{2iu}\,}\, e^{itu}du.
\eeqn
Throughout the paper  $\lg z$ and $\sqrt z$ denote the principal branches in $-\pi <\arg z <\pi$  of the logarithm and the square root, respectively.

 Although one can derive a leading term of $m(t;1)$  directly from (\ref{BE0})    as in \cite{Ubh} (the proofs of Lemmas 4 and 5), here we use the formula
 \beqn\label{2.24} 
 \int_0^\infty N(t)e^{-\la t} dt = \frac{1}{\la-1} - \frac1{\la\lg \la}~~~~~~(\la>0)
\eeqn
(cf. \cite{Ha}, page 196; also \cite{E2}), which somewhat  simplifies the proof.
Put
$$\fa(z)=  \frac{  K_1(\sqrt{2z})}{K_0( \sqrt{2z}) \sqrt{2z}} - \frac{1}{\k}\bigg[  \frac{1}{\k^{-1}z-1} - \frac1{\k^{-1}z \lg(\k^{-1} z)}\bigg].$$
Then (\ref{2.23}) is written as
$$ \int_0^\infty m(t; 1)e^{-\la t}dt = 2\pi   \int_0^\infty N(\k t)e^{-\la  t}dt + 2\pi\fa(\la).$$
Accordingly  (\ref{BE0}) becomes 
\beqn\label{BE}
m(t, 1) - 2\pi N(\k t) = \int_{-\infty}^{\infty}\fa(iu)e^{itu}du.
\eeqn
Here (and in  below) the trigonometric integral  at infinity  is improper.
Note that $\fa(z)$ is analytic on the slit domain $-\pi <\arg z  <\pi$  since $K_0(\sqrt z)$ has no zeros in it: the apparent singularity at $z= \k$ (i.e., that  at $\la=1$ on the right-hand side of (\ref{2.24})) is removable;  also $\fa(iu) $ tends to zero as $u\to\infty$ and is bounded about the origin as will be  observed shortly. 

The idea of the proof of the first formula  in Theorem \ref{thm1} would now be  obvious. The asymptotic behavior of the Fourier integral on the right side of (\ref{BE})  for large values of $t$ depends on that of $\fa(iu)$ near zero, provided that $\fa(iu)$ behaves sufficiently regularly.

In view of  the asymptotic formula 
\beqn\label{K0}
K_\nu(\sqrt{2z})=\bigg(\frac{\,\pi \,}{\, 2\sqrt{2z}\, }\bigg)^{1/2}e^{-\sqrt{2z}}\Bigg(1+\frac{\,\nu^2-\frac14 \, }{2 \sqrt{2z}} + O(|z|^{-1}) \Bigg)~~~~~~~~\mbox{as}~~~|z|\to\infty
\eeqn
$(-\pi<\arg z <\pi,~ \nu\geq 0)$ (cf. \cite{L} (5.11.9)), 
we have 
\beqn\label{f9}
\fa(iu)=\frac{1}{\sqrt{2iu}} - \frac{3}{4i u} + \frac{1}{iu\lg ( iu /\k) } +R(u),
\eeqn
where $R(u) =O(|u|^{-3/2}) $ with its derivatives $R^{(j)}(u)= O(|u|^{-3/2-j})$ as $|u|\to \infty$. 

For the first formula of Theorem \ref{thm1}
it  suffices to show that as $t\to \infty$
\beqn\label{m_0}
\int_{-\infty}^\infty \fa(iu)w(u)e^{itu}du=\frac{4\pi}{t(\lg t)^3}(1+o(1)),
\eeqn
where $w(u)$ is a smooth function that equals 1 in a neighborhood of the origin  and vanishes outside a finite interval, for the $1-w(u)$ part contributes to
the integral at most $O(1/t^N)$ for any $N>1$ so that it  is negligible. 

 Put 
\beqn\label{g(u)}
g(z)=-\lg\, ({\textstyle \frac12}  \sqrt{2z}\,)-\ga = -2^{-1}\lg (\k^{-1} z).
\eeqn 
  By definition
\beqn\label{998}
K_0(z)=\sum_{k=0}^\infty \frac{(z/2)^{2k}}{(k!)^2}\bigg(\sum_{m=1}^k\frac1{m}-\ga-\lg( {\textstyle \frac12} z)\bigg)
\eeqn
and 
$$K_1(z)=\frac1{z} +\frac{z}2  \sum_{k=0}^\infty \frac{(z/2)^{2k}}{k!(k+1)!}\bigg(\lg\frac{z}2 +\ga+ \frac1{2(k+1)}- \sum_{m=1}^k\frac1{m} \bigg)$$
 for $-\pi<\arg z<\pi$.
It follows that for $-1/2<u<1/2$
\beqn\label{g_u0}
K_0(\sqrt{2iu}\,)=g(iu)+2^{-1}iu(g(iu)+1)+R_1(u);
\eeqn
\beqn\label{g_u1}
\frac{K_1(\sqrt{2iu}\,)}{\sqrt{2iu}}=\frac1{2iu} +\frac12 \bigg(-g(iu)+\frac12\bigg)   +R_2(u);
\eeqn
and
\beqn\label{g_u}
\frac1{K_0(\sqrt{2iu}\,)}= \frac1{g(iu)}+\frac{-iu}{2g(iu)}\bigg(1+\frac{1}{g(iu)}\bigg)+ R_3(u).
\eeqn
Here $R_1(u)= O(u^2)\times g(iu)$,  $R_2(u)= O(u)\times g(iu)$ and $R_3(u)={O(u^2)}/{g(iu)}$.  
Substitution of these together with  an easy computation  yields  
 $$\fa(iu) = -\frac12 +\frac{1}{\k} - \frac1{[2g(iu)]^2} +R_4(u)~~~~~~~~~(-1/2<u<1/2)$$
where $R_4(u)=O(u)$ with  the derivatives $R'_4(u)=O(1), R^{(1+j)}_4(u) = O(1/u^j\lg|u|)$ ($j=1,2$).
For evaluation of the contribution of $R_4$ to the Fourier  integral in (\ref{m_0}) we split its range at $u=\pm 1/t$. The inner part plainly gives  the bound $O(1/t^2)$, whereas for the outer part  $\int_{|u|>1/t} R_4(u)w(u)e^{itu}du$ we repeat the integration by parts four times, which leads to the same bound of $O(1/t^2)$ (see  Lemma 2.2 of \cite{U2}  for more details).   
 Hence
$$\int_{-\infty}^\infty \fa(iu)w(u)e^{itu}du = \int_{-\infty}^\infty \frac{-e^{itu}}{[2g(iu)]^2}du + O\bigg(\frac1{t^2}\bigg).$$
 We rewrite the integral on the right-hand  side as 
 $$\frac{\k}{i}  \int_{-i\infty}^{i\infty} \frac{\, -e^{\k tz} \,}{[\lg z]^2}dz$$
  and apply  the Cauchy integral theorem. The integral  along  the  lower (resp. upper) half of the imaginary axis  equals  the one  along the lower (resp. upper) side of the negative real axis in the positive (resp. negative) direction.  Noting $\lg (-x\pm i0)= \lg x \pm i\pi$ for $ x>0$ we then find  the foregoing integral
equal to
$${\k} \int_0^\infty\frac{\, -4\pi (\lg x)e^{-\k t x}\,}{[(\lg x)^2+\pi^2]^2}dx  =\frac{4\pi}{t(\lg t)^3}(1+o(1))$$
(as $t\to\infty$). Thus the relation (\ref{m_0}),  and hence the first formula of Theorem \ref{thm1},  has been proved.
\v2

Consider the case $t\downarrow 0$.  Here we go back to the original inversion integral in  (\ref{BE0}).  From the asymptotic formula (\ref{K0}) we observe as before that the integrand of it involves the singular components $1/\sqrt{2iu}$ and $1/4iu$ that are not integrable at  infinity and we evaluate the contributions of them separately.  To this end 
we bring in a  function $\psi(z)$  by
$$\psi(z)=  \frac{  K_1(\sqrt{2z})}{K_0( \sqrt{2z}) \sqrt{2z}}  -\frac1{\sqrt{2z} }-\frac1{4(z+1)},$$
so that  $\psi(z)= O\Big(|z|^{-3/2}\Big)$ as $z\to \infty$ in the domain $-\pi<\arg z <\pi$ where $\psi(z)$ is regular. 
In addition we have   $\psi(z) =o(1/z)$  if  $z\to 0$ and apply the Cauchy integral theorem  to see that
  the principal value integral $p.v. \int_{-\infty}^\infty \psi(iu) du$ vanishes, so that as $t\downarrow 0$
\beqn \label{t00}
p.v.\int_{-\infty}^\infty \psi(iu)e^{itu}du  = \int_{-\infty}^\infty \psi(iu)\Big[e^{itu}-1\Big] du = O(\sqrt t\,),
\eeqn
where the last equality may be verified by splitting the integral at $u=\pm 1/t$.
  
On the other hand
\beqn\label{t01}
\int_{-\infty}^\infty \frac{1}{\sqrt{2iu}}e^{itu}du =\int_0^\infty\frac{\cos u +\sin u}{\sqrt{u}} du =\sqrt{\frac{2\pi}{t}},
\eeqn
and 
\beqn\label{t02}
\int_{-\infty}^\infty \frac1{4(iu+1)}e^{itu}du = \frac12 \int_{0}^\infty \frac{\cos tu +u\sin tu}{u^2+1} du = \frac{\pi}{2}e^{-t}.
\eeqn
From (\ref{t00}), (\ref{t01}) and (\ref{t02}) we find the formula  as $t\downarrow 0$ as asserted in Theorem \ref{thm1}.

The proof of Theorem \ref{thm1} is complete.

\section{Preliminary results on the first hitting time to a disc}

Here we give several  results  which prepare for  estimation of  the expected area of the sausages for the Brownian bridges. 
 Recall that  $q(\z,t;r) $ denotes the density of the distribution of $\sigma^{(r)}$  the first hitting time at the   disc $U(r)$  of  Brownian motion $B_t$ started at $\z$.
 The following result is proved in \cite{Ubh}.
  \vskip2mm\n
 \begin{Thm}\label{thm0}~   Uniformly for $|\z|>r$,  as $t\to\infty$
\begin{eqnarray}\label{q0}
q(\z,t;r) =   \frac{\lg( \frac12 \k (\z/r)^2)\,}{(\lg (\k t/r^2))^2t} \,e^{- \z^2/2t} + \left\{\begin{array}{ll}   {\displaystyle \frac{2\ga \lg({t}/{\z^2})}{t(\lg t)^3}+ O\bigg(\frac1{t(\lg t)^3} \bigg) }~~~~~ & \mbox{for}~~~\z^2<t, \\[6mm]
{\displaystyle  O\bigg(\frac{1+[\lg (\z^2/t)]^2}{\z^2(\lg t)^3}\, \bigg) }~~~~~& \mbox{for}~~~\z^2 \geq t.
\end{array}\right.
\end{eqnarray}
\end{Thm}

\vskip2mm
Let  $\sigma=\sigma^{(1)}$  and $q(\z,t) =q(\z,t;1)$.  We shall use the estimate of this theorem  in the following slightly reduced form.

\begin{Cor}\label{lem3.1}~
Uniformly for $|\z| >1$, as $t\to\infty$
\beq
q(\z,t)    &=&   \frac{\lg( \frac12 \k (\z/r)^2)\,}{(\lg (\k t/r^2))^2}2\pi p_t(\z)  + O\bigg( \frac{|\lg[(t/\z^2)\vee2]}{t( \lg  t)^3} \bigg)~~~~\mbox{for}~~~~ |\z|^2< 4t\lg (\lg t), \\
&=&  O\bigg( \frac{1}{t( \lg  t)^3} \bigg)~~~~~~~~~~~~~~~~~~~~~~~~~~~~~~~~~~~~~~~~~~\mbox{for}~~~~|\z|^2\geq 4t\lg (\lg t).
\eeq
\end{Cor}
\v2\n
\pf ~ Using the inequality $\lg(\z^2/t) < \z^2/t$  for $\z^2 > t$  the assertion  is immediate from Theorem \ref{thm0}.  \qed

The following crude  bound is often useful for dealing with the case $|\z| > \sqrt{6t\lg \lg t}$. 
\begin{Lem}\label{lem3.11} ~There  is a  constant $c>0$ such that for all $|\z|>1$ and  $t>1$,
 \beqn\label{density_bound} 
q(\z,t) \leq  c\, p_{t+1}(\z). 
\eeqn 
\end{Lem}
\pf ~ For any unit vector $\xi$,  $p_s(\xi)$ depends only on $s$ and we see
$$p_t(\z) =\int_0^t q(\z, t-s)p_s(\xi)ds.$$ 
Hence
$$p_{t+1}(\z) \geq \int_0^1 q(\z,t+1-s)\frac{e^{-1/2s}}{2\pi s}ds$$
and,  applying the parabolic Harnack inequality  (cf. \cite{Mo}) which implies that  $q(\z,t) \leq Cq(\z, t+1-s)$ for $0\leq s\leq 1, |\z|>2$, we obtain the inequality of the lemma. (Note that for $|\z|<2$, a better estimate is given in Corollary \ref{lem3.1}.) \qed

{\sc Remark 3.} ~ On the right-hand side of (\ref{density_bound}) we can replace $p_{t+1}(\z)$ by $p_t(\z)$ for  $|\z| <Mt$, provided that  at the same time $c$ is replaced by a constant $c_M$ which  may depend on $M$. It is warned that the constant $c_M$ actually depends on $M$ and in fact  the ratio $q(\z,t)/p_t(\z)$ tends to infinity whenever $|\z|/t \to\infty$, $t\to\infty$.
\v2

\begin{Lem}\label{lem3.12} ~For all $|\z|>1$ and  $t>0$,
\beqn\label{density bound1}
P_\z[\sigma <t]\leq \sqrt{\frac{2e}\pi} \,\Bigg(\frac{\sqrt t}{|\z|-1} \wedge 1\Bigg)e^{- (|\z|-1)^2/2t}.
\eeqn 
\end{Lem}
\pf ~
The  relation follows
 from the one dimensional result that if $B_t^{(2)}$ denotes the vertical component of $B_t$, 
$ P_{(0,|\z|)}[ B_s^{(2)} > 1~ \mbox{for}~  0<s<t] =2(2\pi t)^{-1/2} \int^\infty_{(|\z|-1)} e^{-u^2/2t}du$. \qed

\begin{Lem}\label{lem3.2}~ There exists $c_1 >0$ such that  for $|\z|>2$ and $t\geq 1$,
$$0 \leq  \frac{d}{dt}E_\z[B_\sigma\cdot \z; \sigma<t] \leq  c_1 \frac{\z^2}{t}p_{t+1}(\z).$$
\end{Lem}
\pf~   
The expectation $E_\z[B_\sigma\cdot \z;  \sigma<t]$ is a radial function of $\z$ and we may suppose $\z$ is on the upper vertical  axis. 
Let $\tau_0$ be the first exit time from $\H\setminus U(1)$, where $\H$ denotes the upper half plane. Let $B_t^{(1)}$ and $B_t^{(2)}$ be  the horizontal  and vertical components, respectively, of $B_t$.  Then by symmetry,  for $\z=(0,y), y>1$,
$$\frac{d}{dt}E_{(0,y)}[B_\sigma\cdot \z; \sigma <t] = y\frac{d}{dt} E_{(0,y)}[ B^{(2)}_{\tau_0}; B_{\tau_0}\in \H, \tau_0  <t].$$ 
The derivative on  the right-hand  side is expressed as   the  integral of the vertical component of  $\xi  \in \H\cap  \partial U(1)$ 
by  $P_{(0,y)}[B_{\tau_0}\in d\xi, \tau_0  \in dt]/dt$.
It therefore suffices to show
that 
\beqn\label{D}
\frac{d}{dt}P_{(0,y)}[B_{\tau_0}\in \H, \tau_0  <t]  \leq C \frac{y}{t^2}e^{-y^2/2t}.
\eeqn
For  verification  let $D = \H-(0,1) =\{(x,y-1): (x,y) \in {\bf H}\}$. Then, on denoting  by $\tau_D=\tau(D)$  the first exit time from $D$,
\beq
&&P_{(0,y)}[\tau_{D}\in dt, \, -1<  B^{(1)}_{\tau(D)} <1] \\
&&\geq \int_{|\xi|=1, \xi\in \H} \int_0^t ds  \frac{d}{ds}P_{(0,y)}[B_{\tau_0}\in d\xi, \tau_0\leq s]P_\xi[ B^{(1)}_{\tau(D)}\in (-1,1), \tau_D+s\in dt]. 
\eeq
Since
$$P_{(x,y)}[\tau_{D}\in dt,  B^{(1)}_{\tau(D)}\in  du ]= p_t((x-u,y+1))\frac{y+1}{t}dtdu,$$
 we obtain
\beq
2p_t((0, y+1)) \frac{y+1}{t}  &\geq& \int_{|\xi|=1, \xi\in \H} \int_0^t   \frac{d}{ds}P_{(0,y)}[B_{\tau_0}\in d\xi, \tau_0\leq s] \frac{2p_{t-s}(\xi+(0,1))}{t-s} ds\\
&\geq&  c \int_{t-1/2}^{t-1/4}   \frac{d}{ds}P_{(0,y)}[B_{\tau_0}\in \H, \tau_0\leq s]  ds.
\eeq
for some universal constant $c>0$. Thus the  bound of   (\ref{D})   follows in view of the parabolic Harnack inequality as before.

\section{Wiener sausages for Brownian bridges}
For $\x\in\R^2$, put
\beq
F(t, \x ; r) = \int_{|\z|\geq r}d\z \int_0^t \int_{|\xi|=r} P_\z[\sigma_{r} \in ds, B_{\sigma_{r}} \in d\xi ]p_{t-s}(\x-\z-\xi).
\eeq
  Then  
\beqn\label{Area}
E[{\rm Area} (S_t^{(r)})\,|\, B_t=\x]  =\frac1{p_t(\x) }  F(t, \x; r)+  \pi  r^2.
\eeqn
From the scaling property of Brownian motion it follows that
$$F(t,\x; r) = F(t/r^2,\x/r; 1)$$
and we have only to  consider  the case  $r =1$ as in the proof of Theorem \ref{thm1}.
Put
$$F_0(t, \x) = \int_{|\z|\geq 1} d\z\int_0^tP_\z[\sigma \in ds]p_{t-s}(\z-\x) =\int_{|\z|\geq 1}d\z \int_0^t q(\z,s) p_{t-s}(\z- \x)ds.
$$
We first consider the case $\x =0$ in Subsection 4.1. The general case, dealt with in  Subsection 4.2, is based on the proof of this special case but need additional estimations of a certain integral  that are 
 involved and partly delicate.
\v2\n
\subsection{ The case $\x=0$}
\begin{Lem}\label{lem3.3}~ Let $F$ and $F_0$ be as above. Then   $F(t,0;1)= F_0(t,0) + O(1/t)$.
\end{Lem}
\pf~ Let $a$ be a constant larger than 1. In this proof it may be arbitrarily fixed (eg. $a=2$; in the case $\x\neq 0$ treated later we shall take $a=\x^2$).  We split  the range of integration by the surfaces $s=t-a$ and $|\z| = \sqrt{4(t-s)\lg (t-s)}$  put 
\beq
&& D_0=D_0(a)= \{(s,\z): t-a <s  \leq t, |\z|>1\} ~~~\\
&&D_> =D_>(a) =\{(s,\z): 0<s\leq t-a,  |\z| > \sqrt{4(t-s)\lg (t-s)}, ~ |\z|>1\},  \\
&&D_< =D_<(a)= \{(s,\z): 0<s\leq t-a, |\z| \leq \sqrt{4(t-s)\lg (t-s)}, ~ |\z|>1\}.
\eeq
Accordingly we break the rest of the proof into three parts. Some estimates obtained below are more accurate than necessary for the poof of the 
present lemma, but they are needed  in the proof of  Lemma \ref{lem3.5} essential for our proof of Theorem \ref{thm2}. 

\v2
{\it Part 1} :$D_0$.~
The contribution to the integrals defining $F$ and $F_0$ from $D_0$ is $O(1/t(\lg t)^2)$. Indeed on the one hand we can  use the bound (\ref{density_bound})
for the  integration  w.r.t.  $\z$ over the range $|\z| >\sqrt{4t\lg \lg t}$, in which $p_{t+1}(\z) = O(1/t(\lg t)^2)$.   On the other hand, on applying   Corollary \ref{lem3.1}  the  contribution of the other part is at most a constant multiple of 
  \begin{eqnarray}\label{K}
 && \int_{t-a}^t   \frac{ds}{s (\lg s)^2} \int_{1< |\z| < \sqrt{4t\lg \lg t} }\,(\lg |\z|) \sup_{|\xi|=1} p_{t-s}(\z-\xi)d \z\nonumber\\
 &&\leq \frac{C}{t(\lg t)^2}.
  \end{eqnarray}

{\it Part 2} :$D_>$. ~  For $(s,\z)\in D_>$,  according to  Corollary \ref{lem3.1}  $q(\z,s)\leq C/s(\lg s)^3$  if $a<s<t/2$ and   $q(\z,s)\leq C/s\lg s$  if $s>t/2$. Hence both the contributions of $D_>$ to $F$ and $F_0$ are  dominated from above by a constant multiple of 
  \beq
 && \int_4^{t-a} \Bigg[\frac{\1(s\leq t/2)}{s(\lg s)^3} + \frac{\1(s>t/2)}{s\lg s}\Bigg] ds\int_{|\z|>\sqrt{4(t-s)\lg (t-s)}} \sup_{|\xi| =1}p_{t-s}(\z-\xi)d\z \\
 &&\leq \frac{C'}{t^2} + C' \int_{t/2}^{t-a} \frac{ds}{(s\lg s)(t-s)^2}\\ 
 &&= O\bigg(\frac1{t\lg t}\bigg).
  \eeq
 On using  (\ref{density bound1}) the integral on $0<s\leq 4$ is readily evaluated  to be at most $O(1/t^2)$.
\v2
{\it Part 3} :$D_<$. ~  Observe that if  $(s,\z) \in D_<$, then
\begin{eqnarray}\label{B}
p_{t-s}(\z-\xi) - p_{t-s}(\z) &=& p_{t-s}(\z)( e^{[\z\cdot \xi -\frac12 ]/(t-s)}-1) \nonumber \\
&=& p_{t-s}(\z)\frac{\z \cdot \xi -\frac12}{t-s} + \eta_t(s,\z,\xi) 
\end{eqnarray}
with 
$$\eta_t(s,\z,\xi)= p_{t-s}(\z)\times O\bigg(\frac{\z^2 +1}{(t-s)^2} \bigg).$$

The  integral for the difference $F(t,0;1)-F_0(t,0) $ restricted to
$$D_<^{0}:=\{(s,\z) \in D_<: 0<s\leq 4\}$$
 is at most $O(1/t^2)$ in view of  (\ref{density bound1}).
Thus we may restrict the integral to the range  $D_< \setminus D_<^0$.

Consider the contribution of  the first term on the right side of (\ref{B}).  Applying Lemma \ref{lem3.2} with the help of Corollary \ref{lem3.1} (for $|\z|< 2$)  
we deduce that
\begin{eqnarray}\label{A1}
&&\int_{D_< \setminus D_<^0}ds  d\z  \Bigg| \int_{|\xi|=1}   p_{t-s}(\z) \frac{\z \cdot \xi}{t-s}\cdot \frac{P_\z[\sigma \in ds, B_{\sigma} \in d\xi ] }{ds} \Bigg|  \nonumber\\
&&\leq  C \int_4^{t-a} ds\int_{|\z|>2}  \frac{p_{t-s}(\z)\z^2 p_{s+1}(\z)}{(t-s)s}d\z + C\int_4^{t-a}  \frac{ds}{(t-s)^2s(\lg s)^2} \leq
\frac{C'}{t}.
\end{eqnarray}
Here we have applied the trivial bound $p_{t-s}(\z)< 1/(t-s)$ (for the integral on $s<t/2$) as well as  $p_{s+1}(\z) \leq1/s$ (for that on $s>t/2$).  The part involving  the term $1/2$ is  evaluated to be $O(1/t)$  by dominating  $q(\z, s) $ by $Cp_{s+1}(\z)$ for $s<t/2$ and by  $C/s\lg s$ for $s>t/2$.  

For the proof of the lemma it now  suffices to show 
that
\beqn\label{D1}
\int_{D_< \setminus D_<^0}ds  d\z   \int_{|\xi|=1}   \eta_t(s,\z,\xi)    \cdot \frac{P_\z[\sigma\in ds, B_{\sigma} \in d\xi ] \,}{ds}
 = O\bigg(\frac1{t\lg t}\bigg).
\eeqn

To this end  we further split the region $D_<$  by the surfaces $s=t/2$ and  $ |\z|= 1+\sqrt{4s\lg\lg s}\,$ ($4<s<t/2$).  Denote by  $J^{{\rm RH}}(a)$, $J_{{\rm LH}}^{(>)}$, $J_{{\rm LH}}^{(<)}$   the double integrals  on  the  regions   
\beq  
&&D_<^{{\rm RH}}(a):=  \{(s,\z) \in D_<: t/2\leq s<t-a\},\\  
&&D_{{\rm LH}}^{(>)} :=\{(s,\z) \in D_<: 4<s< t/2,  |\z|\geq 1+\sqrt{4s\lg\lg s}\,\},\\
&&D_{{\rm LH}}^{(<)} :=\{(s,\z) \in D_<: 4<s< t/2,  |\z| <1+\sqrt{4s\lg\lg s}\,\},
\eeq
  respectively. 
Employing Corollary  \ref{lem3.1}   we  obtain 
\begin{eqnarray}\label{A0}
 J^{{\rm RH}}(a) & < & C\int _{t/2 }^{t-a} \frac{ds}{s\lg s} \int_{|\z|< \sqrt{4(t-s)\lg (t-s)}} \, \frac{\z^2}{(t-s)^2}\, p_{t-s}(\z)d\z \nonumber\\
 &<&\frac{C'}{t \lg t}  \int _{t/2 }^{t-a} \frac{ds}{t- s} =  O\bigg(\frac1{\,t\,}\bigg),
\end{eqnarray}
and, similarly but dominating $p_{t-s}(\z)$ simply by $1/ (t-s)$,
\begin{eqnarray}\label{A}
J_{{\rm LH}}^{(<)}& < & C\int _{4}^{t/2 }\frac{ds}{(t-s)^3} \int_{|\z|< 1+\sqrt{4s\lg\lg s}}\,  \frac{\z^2}{\lg s} p_s(\z)d\z \nonumber\\
 &<&\frac{C'}{t^3}\int_4^{t/2} \frac{s}{\lg s}ds =  O\bigg(\frac1{t\lg t}\bigg).
\end{eqnarray}
 With the help of (\ref{density_bound}) and  (\ref{q0}) as well as  (\ref{B})  we deduce that
\beq
J_{{\rm LH}}^{(>)} &<&  C \int_4^{t/2}ds\Bigg[ \int_{|\z|>1+\sqrt{6s\lg\lg s} }\, p_{s+1}(\z) \frac{|\z|^2 p_{t-s}(\z)}{(t-s)^2}d\z\\
&&~~~~~~~~~~~~~~~+  \int_{|\z|>1+\sqrt{4s\lg\lg s} }\,\bigg\{\frac{p_s(\z)}{\lg s} +  \frac{(\lg\lg s)^2}{\z^2(\lg s)^3}\bigg\} \frac{\z^2p_{t-s}(\z)}{(t-s)^2}d\z  
 \Bigg] \\
& <&  C'\int_4^{t/2} ds  \Bigg[ \frac{1}{(\lg s)^2t^2} +  \int_{|\z|>1+\sqrt{4s \lg\lg s} }\, \frac{\z^2p_{s}(\z) d\z}{t^3\lg s} + \int_{|\z|>1+\sqrt{6s\lg\lg s} }\, \frac{\z^2p_{s}(\z) d\z}{t^3}  \Bigg]\\
&\leq& 
  \frac{C'' }{t(\lg t)^2}.  
\eeq
Thus (\ref{D1}) has been proved.
The proof of the lemma is complete.   \qed
\v2

\begin{Lem} \label{lem3.4} ~  $$F_0(t,0) = N(\k t)  -\pi p_t(1) + O(1/t(\lg t)^2).$$
\end{Lem}
\pf ~ 
The  Laplace transform of $F_0$ is given by
\beq
\int_0^\infty F_0(t,0)e^{-\la t} dt&=&\int_0^\infty e^{-\la t}dt\int_{|\z|\geq 1} d\z\int_0^tP_\z[\sigma \in ds]p_{t-s}(\z)\\
&=&\frac1{\pi} \int_{|\z|\geq 1} \Big[K_0(|\z|\sqrt{2\la})\Big]^2d\z\,  \frac{1}{K_0(\sqrt{2\la})}
\eeq
(cf.  \cite{IM} Sect.7.2).  Fortunately  we have the identity 
$$rK_0^2(ra) =\frac12\frac{d}{dr}\Big[ r^2\Big(K_0^2(ra) - K_1^2(ra)\Big)\Big]$$
for any constant $a>0$,  so that
$$\int_0^\infty F_0(t,0)e^{-\la t}dt  = -K_0(\sqrt{2\la}) +  \frac{[K_1(\sqrt{2\la})]^2}{K_0(\sqrt{2\la})}.$$
Hence
$$F_0(t,0)= - \pi p_t(1) + \frac{1}{2\pi }\int_{-\infty}^{\infty}\frac{[K_1(\sqrt{2iu})]^2}{ K_0(\sqrt{2iu})}e^{itu}du.$$
For evaluation of the last integral we can proceed as in the preceding section. Put
$$E(z)=  \frac{[K_1(\sqrt{2z})]^2}{K_0(\sqrt{2z})} - \frac{1}{\k}\bigg[  \frac{1}{\k^{-1}z-1} - \frac1{\k^{-1}z \lg(\k^{-1} z)}\bigg].$$
Then for $-1/2<u<1/2$, $E(iu) = -1 +\frac1{\k} +[4g(iu)]^{-1} - [2g(iu)]^{-2} + R(u)$
with $R(u)=O(u)$ and as before (see (\ref{2.24}) and the ensuing discussion up to (\ref{BE}))
we have
$$F_0(t,0) - N(\k t) +\pi p_t(1) =\frac1{2\pi}\int_{-\infty}^\infty E(iu)e^{itu}du = O\bigg(\frac1{t(\lg t)^2}\bigg).$$
The proof of the lemma is complete. \qed

Combining  Lemmas \ref{lem3.3} and \ref{lem3.4} with (\ref{Area}) leads to the assertion of Theorem \ref{thm2}  with $\x=0$.
\v2\n
\subsection{ The case $\x\neq 0$.}
We first extend Lemma \ref{lem3.3}.
\begin{Lem} \label{lem3.50} ~  For each $M\geq 1$ uniformly for $|\x|< M\sqrt{t}$, as $t\to \infty$  
$$F(t,\x;1) = F_0(t, \x) + O(1/{t}).$$
\end{Lem}
\pf~ First of all we verify that 
\beqn\label{E}
\int_{|\x-\z|<1, |\z|>1}d\z\int_{t-1}^t P_\z[\sigma\in ds, B_\sigma \in d\xi]p_{t-s}(\z-\x-\xi)\leq \frac{C}{t\lg t}.
\eeqn
For verification we use the inequality 
\beq
P_\z[\sigma\in ds, B_\sigma \in d\xi]\leq \int d {\bf y} p_1({\bf y}-\z) P_{{\bf y}}[\sigma +1 \in ds, B_\sigma \in d\xi].
\eeq
If  $|\x-\z|<1$,  $p_1({\bf y}-\z)$ is dominated by $e^{-({\bf y}-\x)^2/4}$.  Hence,  performing the  integration by $\z$ first, we see that the repeated integral in (\ref{E}) is at most   $\int_{t-1}^t ds \int  e^{-({\bf y}-\x)^2/4}q({\bf y}, s-1)d{\bf y}$, which is evaluated to be $O(1/t\lg t)$ according to  Corollary \ref{lem3.1}.

The rest of the proof proceeds in parallel with  that of Lemma \ref{lem3.3}.  Of course we must   replace $p_{t-s}(\z-\xi)$ and  
$p_{t-s}(\z)$  by $p_{t-s}(\z-\x-\xi)$ and by  $p_{t-s}(\z-\x)$, respectively.  We accordingly  modify  the definition of  $D_<$  as
$$D_< = \{(s,\z): 0<s <t-a, |\z-\x| \leq \sqrt{4(t-s)\lg (t-s)}\}$$
 and $D_>$
  analogously. Then Parts 1 and 2 are treated without any change except that for Part 1 we use  (\ref{E}).
For Part 3 we consider
\begin{eqnarray}\label{B1}
p_{t-s}(\z- \x-\xi) - p_{t-s}(\z-\x) &=& p_{t-s}(\z-\x)( e^{[(\z-\x)\cdot \xi -\frac12 ]/(t-s)}-1) \nonumber \\
&=& p_{t- s}(\z-\x) \frac{(\z-\x) \cdot \xi -\frac12}{t-s} + \eta_{t,\x}(s,\z,\xi) 
\end{eqnarray}
with 
$$\eta_{t,\x}(s,\z,\xi)=  p_{t- s}(\z-\x) \times  O\bigg(\frac{|\z-\x|^2 +1}{(t-s)^2} \bigg)$$
in place of (\ref{B}).  In the first term on the last member of  (\ref{B1})  $\x\cdot\xi$ may be replaced by $(\x\cdot \z)( \z\cdot \xi)/\z^2$
since by symmetry  the component of $\xi$  perpendicular to $\z$ vanishes after integration. It follows  that $(\z- \x)\cdot\xi$ can be replaced by $(\z\cdot\xi)\z\cdot (\z- \x)/|\z|^2$.   
 Keeping this in mind and  employing (\ref{density bound1}) together with Lemma \ref{lem3.2} and Corollary \ref{lem3.1},  we deduce that 
\begin{eqnarray}\label{A10}
&&\int_{D_<}ds  d\z  \Bigg| \int_{|\xi|=1}   p_{t-s}(\z-\x) \frac{(\z-\x) \cdot \xi -\frac12}{t-s}\cdot \frac{P_\z[\sigma \in ds, B_{\sigma} \in d\xi ] }{ds} \Bigg|  \nonumber\\
&&\leq \frac{C}{t} + \int_{a}^{t-a}ds  \int_{|\z|\geq 2} \frac{|\z\cdot (\z-\x)|}{(t-s)s}\cdot  p_{t-s}(\z-\x) p_s(\z)d\z\nonumber\\
&&~~~~ +
C\int_{|\z|<2}  d\z \int_a^{t-a} \, |\z - \x| \, e^{-|\z-\x|^2/2(t-s)} \frac{ds}{(t-s)^2s(\lg s)^2}  \nonumber\\
&&\leq \frac{C'}{t}.
\end{eqnarray}
Here  the last inequality  is due to the observation   that the first  integral  in the middle member restricted on the interval $s\in [a,t/2]$ is at most
a constant multiple of 
$$\frac1{t^2}\int_a^{t/2} ds\int_{\R^2}  \frac{\,\z^2 + |\z||\x| \,}{s} p_s(\z)d\z \leq \frac{1}{t} + \frac{|\x|}{t^2}\int_0^{t/2}\frac1{\sqrt{s}}\,ds  =O\bigg(\frac1{\,t\,}\bigg)  $$
and similarly for the other interval.

The rest of the proof is similarly dealt with and hence omitted.
\qed

\begin{Lem} \label{lem3.5} ~ Foe each $M\geq 1$ uniformly for $|\x|<M \sqrt{t}$, as $t\to \infty$ 
\beq
F_0(t,\x) - F_0(t,0) &=&  2\pi  \Big[p_{t}(\x) -p_{t}(0)\Big]  \frac{t}{\lg t}  
+ \frac{ \x^2}{2t(\lg t)^2}\bigg[\lg\bigg(\frac{t}{\x^2\vee1}\bigg) +O(1)\bigg]  \\
&&+ ~ O\bigg(\frac{1}{t\lg t} \bigg).
\eeq
\end{Lem}
The right-hand side of the formula of this lemma  is suggested  by that of Theorem \ref{thm2} and the latter in turn is  by the corresponding one   in \cite{U1} for  random walks. The course of  proof given below  is often  steered   by the asserted formula.
\v2\n
{\it Proof of Lemma \ref{lem3.5}.} ~We suppose $|\x|<  \sqrt{t}\,/2$ for simplicity of the description of the proof (otherwise one may take $a= \x^2/2M^2$ instead of $a=\x^2$ in below). Suppose also  $|\x| \geq 2$,
 the  same proof of Lemma \ref{lem3.3} being  applicable  with little alteration if $|\x|<2$ .   We make computations analogous to those  in the proof of Lemma  \ref{lem3.3}.  We substitute  $a=\x^2$ in the definitions of $D_0, D_>$ and $D_<$ 
given therein and denote the resulting regions by $D_0(\x^2)$, $D_>(\x^2)$ and $D_<(\x^2)$, respectively; denote the contributions from them to  the difference $F_0(t,\x)-F_0(t,0)$ by $I[D_0(\x^2)], I[D_>(\x^2)]$ and $I[D_<(\x^2)]$, respectively.  

For the estimate   given in Part 2 of the proof of Lemma \ref{lem3.3}   where the integral over $D_>(a)$ is dealt with we have the same bound $O(1/t\lg t)$, namely  
 \beqn\label{E1}
 I[D_>(\x^2)] =O(1/t\lg t),
 \eeqn
  since $|\z-\x| \geq |\z|- |\x| \geq  \sqrt{ 3 t\lg t}$ if $t/2 < s <t-\x^2$ and $(s,\z)\in D_>(\x^2)$ (the integral on $0<s<t/2$ is dealt with in the same way as before).
  
  The rest of the proof  is somewhat involved and
  broken  into three parts {\bf 1} to {\bf 3}. 
 \v2
{ \bf 1.} {\it Estimation of}  $I[D_0(\x^2)]$. 
  There exists a constant $C$ such that for  $|\x| >1$,
 \beqn\label{I}
 \int_0^{\x^2}du \int_{\R^2} \Big| \lg \z^2- \lg |\z- \x|^2 \Big|
p_u(\z)d\z < C\x^2
\eeqn
as is proved shortly. 
For the integral that gives  $I[D_0(\x^2)]$   we may restrict $(s,\z)$  to $\{4<t-s< \x^2,  |\z| <\sqrt{4s\lg\lg s}\}$, the contribution of the remainder being at most $O(\x^2/t(\lg t)^2)$ as discussed before. In view of this   as well as of  (\ref{I}) we may write it as
\beq
&&\int_{t-\x^2}^{t-4}ds\int_{|\z| <\sqrt{4 s\lg\lg s} } \frac{2\pi p_s(\z)}{(\lg s)^2}\Big[ (\lg |\z- \x|^2)p_{t-s}(\z-\x) - (\lg \z^2) p_{t-s}(\z)\Big]d\z \\
&&+ ~ O\bigg(\frac{\x^2}{t(\lg t)^2}\bigg).
\eeq
Substituting from the identities  $\lg |\z-\x|^2 =\lg (t-s) +\lg [|\z-\x|^2/(t-s)]$ and $\lg \z^2=\lg (t-s) +\lg [\z^2/(t-s)]$ and noting that  $\int \lg (\z^2/(t-s)) p_{t-s}(\z)d\z$  equals $\int (\lg \z^2)p_1(\z)d\z$, a finite constant  we find that
\beqn\label{1.00}
I[D_0(\x^2)] = 2\pi \Big[p_{t}(\x) -  p_{t}(0)\Big] \int_{t-\x^2}^{t-4}\frac{\lg(t-s)ds}{(\lg s)^2}+ O\bigg(\frac{\x^2}{t(\lg t)^2}\bigg).
\eeqn
This  is not negligible and to  contribute to the leading term on the right-hand side of  the formula of the lemma.
\v2\n
{\sc Proof of} (\ref{I}).   If the range of integration  is restricted to $|\z| \wedge  |\z-\x| >|\x|/2$, then  $\Big| \lg \z^2- \lg |\z- \x|^2 \Big| $ is uniformly  bounded and the corresponding (repeated)  integral is obviously bounded by a positive multiple of $\x^2$. The integral corresponding to $\{|\x-\z| \leq |\x|/2\}$ can  be readily evaluated to be  $O(\x^2)$.  For the range  $|\z| \leq |\x|/2$ it suffices to show that
\beqn\label{I1}
\frac1{\x^2}\int_0^{\x^2}\frac{du}{u}\int_{|\z|<|\x|} \bigg(\lg \frac{|\x|}{|\z|}\bigg) e^{-|\z|^2/2u}d\z < C'.
\eeqn
Changing the variables of integration according to  $s=u/\x^2$ and  ${\bf y}=  \z/|\x| $ transforms  the left-hand side of (\ref{I1})  to 
$\int_0^1 s^{-1}{ds}\int_{|{\bf y}|< 1}  (\lg |{\bf y}|) e^{-{\bf y}^2/2s} d{\bf y},$ which is independent of $x$.  By another change of  variables this last integral is further transformed to
$$\int_1^\infty \frac{ds}{s^2}\int_{|{\bf y}|<\sqrt s} \bigg(\lg \frac{\sqrt s}{|{\bf y}|}\bigg) e^{-{\bf y}^2/2} d{\bf y},$$
which is certainly finite. Hence we have (\ref{I1}).
 
\v2
 
 {\bf 2.}  {\it Estimation of}  $I[D_<(\x^2)]$.
For $(s,\z)\in D_<(\x^2)$,  instead of  (\ref{B}) we have
  \begin{eqnarray}\label{L}
p_{t-s}(\z-\x) - p_{t-s}(\z) 
&=& p_{t-s}(\z)( e^{[\z\cdot \x -\frac12 \x^2]/(t-s)}-1) \nonumber \\
&=& p_{t-s}(\z)\Bigg[ \frac{\z \cdot \x -\frac12 \x^2}{t-s} +\frac{(\z\cdot \x)^2}{2(t-s)^2} -\frac{(\z\cdot \x)\x^2}{2(t-s)^2} \Bigg] \nonumber \\
&&~~~+ \eta_{t, \x}(s,\z)
\end{eqnarray}
with 
$$\eta_{t,\x}(s,\z)= p_{t-s}(\z) \frac{\x^2}{t-s} \times 
O\bigg(\frac{ \x^2}{t-s} +\frac{ |\z|^4}{(t-s)^2}\bigg).$$
Since the double  integral  of  the terms that are linear in $\z$    inside the big square brackets in (\ref{L})  vanishes by skew symmetry, 
the substantial part  of $p_{t-s}(\z-\x) - p_{t-s}(\z)$ reduces to $p_{t-s}(\z)\times L_{t,\x}(\z)$, where
\beqn\label{F}
L_{t,\x}(s,\z) =  -\frac{\x^2}{2(t-s)}  +\frac{(\z\cdot \x)^2}{2(t-s)^2}+\frac{\x^2}{t-s} \times 
O\bigg(\frac{ \x^2}{t-s} +\frac{ |\z|^4}{(t-s)^2}\bigg).
\eeqn

 Let  $D^0_<$, $D^{{\rm RH}}_<(\x^2)$, $D^{(>)}_{\rm LH}$ and $D_{\rm LH}^{(<)}$ be defined as in Part 3 of the proof of Lemma \ref{lem3.2} and denote by   $I[D^0_<]$, $I[D_<^{{\rm RH}}(\x^2)]$ etc. be corresponding double integrals on them, so that 
 \beqn\label{E2}
 I[D_<(\x^2)] =I[D_<^0]+ I[D_<^{{\rm RH}}(\x^2)]+ I[D_{\rm LH}^{(<)}] + I[D_{\rm LH}^{(>)}].
 \eeqn
  The integrals $I[D_<^0],  I[D_{\rm LH}^{(<)}] $ and $ I[D_{\rm LH}^{(>)}]$ are independent of $\x$ as the notation suggests.
\v2
 {\bf 2.1.} {\it Estimation of} $I[D_<^{{\rm RH}}(\x^2)]$.  By definition 
  \beqn\label{(3.1)}
I[D_<^{{\rm RH}}(\x^2)]= \int _{t/2 }^{t-\x^2} ds \int_{|\z|< \sqrt{4(t-s)\lg (t-s)}} q(\z,s)\, \Big[p_{t-s}(\z-\x) -p_{t-s}(\z)\Big] d\z. 
\eeqn

Let $(s,\z)$ be in the region of integration of this integral.
By what is remarked right above,   $p_{t-s}(\z-\x) -p_{t-s}(\z)$ may be replaced by  $p_{t-s}(\z)L_{t,\x}(s,\z)$, which is bounded in absolute value  by a constant multiple of  
$$p_{t-s}(\z)(t-s)^{-1}\x^2\Big[1 + (t-s)^{-2}|\z|^4\Big].$$
 Keeping  this in mind   
 one  replaces the upper limit $\sqrt{4(t-s)\lg (t-s)}$ of the inner integral in (\ref{(3.1)}) by $\sqrt{8(t-s)\lg\lg (t-s)}$. Since  $q(\z,s)\leq C(\lg |\z|) /s(\lg s)^2$ in view of  Corollary \ref{lem3.1} and the integral of $(\lg |\z|)[1+(t-s)^{-2} \z^4] p_{t-s}(\z)$ w.r.t. $\z$ over ${|\z|>\sqrt{8(t-s)\lg\lg (t-s)}}$ is at most  $1/[\lg (t-s)]^2$, the error given rise to  by  this replacement is bounded above by 
 $$C'\x^2 \int_{t/2}^{t-\x^2} \frac{ds}{s(\lg s)^2  [(t-s)(\lg (t-s))^2]}= O\bigg(\frac{\x^2}{t(\lg t)^2}\bigg).$$
 This together with Corollary \ref{lem3.1} shows that 
     \begin{eqnarray}\label{(3.11)}
I[D_<^{{\rm RH}}(\x^2)]&=& \int _{t/2 }^{t-\x^2} ds \int_{|\z|< \sqrt{8(t-s)\lg\lg (t-s)}}\frac{\lg \z^2}{(\lg s)^2} 2\pi p_s(\z)\, \Big[p_{t-s}(\z-\x) -p_{t-s}(\z)\Big] d\z \nonumber\\
&& +  ~O\bigg(\frac{\x^2}{t(\lg t)^2}\bigg). 
\end{eqnarray}

Now, substituting the decomposition $\lg \z^2 = \lg (t-s) + \lg[\z^2/(t-s)]$, we have
\beqn\label{(3.12)}
I[D_<^{{\rm RH}}(\x^2)]=
2\pi  \Big[p_{t}(\x) -p_{t}(0)\Big]  \int _{t/2 }^{t-\x^2} \frac{\lg (t-s)}{(\lg s)^2} ds  + K(\x,t) + O\bigg(\frac{\x^2}{t(\lg t)^2}\bigg),
\eeqn
where 
$$K(\x,t) =\int _{t/2 }^{t-\x^2} \frac{ds}{(\lg s)^2} \int_{|\z|< \sqrt{8(t-s)\lg\lg (t-s)) }} \Bigg(\lg\frac{\z^2}{ t-s}\Bigg) 2\pi p_s(\z) p_{t-s}(\z)  L_{t,\x}(s,\z)d\z.$$
On writing $u$ for $t-s$, using the identity
$$2\pi p_u(\z)p_{t-u}(\z) =t^{-1} p_{(t-u)u/t}(\z),$$
 noting that one may replace $\lg (\z^2/u)$ by $\lg\Big( \z^2/[u(t-u)/t]\Big)$ in the present estimation, and changing the variables of integration according to ${\bf y}=\z/\sqrt{s(t-s)/t}$,  the inner integral above becomes
\beq
 && \int_{|\z|< \sqrt{8u\lg\lg u}} \, \frac{\lg (\z^2/u)}{t}\, p_{(t-u)u/t}(\z)L_{t,\x}(t-u,\z)   d\z \\
 &&=  \int_{|{\bf y}|< \sqrt{8(t/(t-u))\lg\lg u}} \, \frac{\lg {\bf y}^2}{t}\, \bigg[-\frac{\x^2}{2u\,} + \frac{(t-u)({\bf y}\cdot \x)^2}{2ut}\bigg] p_{1}({\bf y})d{\bf y}  + R \\
 &&=\bigg(-\frac1{2tu} b_1+  \frac{t-u}{4ut^2}b_2\bigg)\x^2 + O\bigg(  \frac{\x^2}{tu(\lg u)^{3/2}} \bigg) +R
 \eeq
(under the restriction $x^2< u<t/2$), 
where  $R$ denotes the contribution of the error term in $L_{t,\x}(t-u,\z)$  and we have used the identity  $\int ({\bf y}\cdot \x)^2f(|{\bf y}|)d{\bf  y} =2^{-1}\x^2\int {\bf  y}^2f(|{\bf y}|) d{\bf y}$ valid for a non-negative function $f(r)$ and
$$b_1 =  \int_{\R^2}(\lg {\bf y}^2)p_1({\bf y})d{\bf y}~~~~\mbox{and}~~~~  b_2= \int_{\R^2}{\bf  y}^2(\lg {\bf y}^2)p_1({\bf y})d{\bf y}.$$
It follows that $b_2-2b_1 = 2$ and $|R|\leq  C|K(\x,t)| \x^2/t$. Noting that $\int_{\x^2}^{t/2} u^{-1}(\lg (t-u))^{-2} du = (\lg t/\x^2)(\lg t)^{-2} +O((\lg t)^{-2})$ we infer 
\begin{eqnarray*}\label{F2}
K(\x,t)
=\frac{ \x^2\lg(t/\x^2)}{2t(\lg t)^2} + O\bigg(\frac{\x^2}{t(\lg t)^2}\bigg),
\end{eqnarray*}
and substitution into (\ref{(3.12)}) yields
\beqn\label{2.10}
I[D_<^{{\rm RH}}(\x^2)] =
2\pi  \Big[p_{t}(\x) -p_{t}(0)\Big]  \int _{t/2 }^{t-\x^2} \frac{\lg (t-s)}{(\lg s)^2} ds  + \frac{ \x^2\lg(t/\x^2)}{2t(\lg t)^2}+ O\bigg(\frac{\x^2}{t(\lg t)^2}\bigg). 
\eeqn
\v2
 {\bf 2.2.} {\it Estimation of} $I[D_{\rm LH}^{(<)}(\x^2)]$.  This part is similar to the preceding one but much simpler. This time we substitute $\lg \z^2 = \lg s + \lg (\z^2/s)$ to see 
 that
 \begin{eqnarray}
 I[D_{\rm LH}^{(<)}(\x^2)] &=&  \int _{4 }^{t/2} ds \int_{|\z|< \sqrt{4s \lg \lg s}} \, q(\z,s)  \Big[p_{t-s}(\z-\x) -p_{t-s}(\z)\Big]\, d\z \nonumber \\
 &=& 2\pi \Big[p_t(\x)-p_t(0)]  \int _{4 }^{t/2} \frac{ds}{\lg s} + R
 \label{2.20}
 \end{eqnarray}
 with 
\beq |R| &\leq&  C\int_{4 }^{t/2} \frac{ds}{(\lg s)^2}\int_{\R^2} p_{t-s}(\z) \Bigg[\frac{\x^2}{t-s} + \frac{\z^2\x^2}{(t-s)^2}\Bigg] \bigg(\lg\frac{\z^2}{s}\bigg) p_s(\z)d\z \\
&\leq& \frac{C' {\x^2}}{{t(\lg t)^2}}.
 \eeq 
\v2
 {\bf 2.3.} {\it Estimation of} $I[D_{\rm LH}^{0}(\x^2)]$ and $I[D_{\rm LH}^{(>)}(\x^2)]$.  
  It is readily checked that  $I[D_<^{0}(\x^2)]  = O(1/t(\lg t)^2)$. By the same argument as made for the estimation of $J_{\rm LH}^{(>)}$ in the proof of Lemma \ref{lem3.3} we deduce that 
 \begin{eqnarray} \label{2.30}
 I[D_{\rm LH}^{(>)}(\x^2)] =  O\bigg(\frac{\x^2}{t(\lg t)^2} \bigg). 
 \end{eqnarray}
 \v2
{\bf 3}. {\it Completion of Proof.} ~Combining (\ref{1.00}), (\ref{2.10}),  (\ref{2.20}), (\ref{2.30}) and (\ref{E2})  we  deduce  that 
 \beq
I[D_<(\x^2)]  &=&  2\pi  \Big[p_{t}(\x) -p_{t}(0)\Big] \Bigg[\int _{t/2 }^{t} \frac{\lg (t-s)}{(\lg s)^2}ds  +\int_0^{t/2}\frac{ds}{\lg s} \Bigg]\\
 && + \frac{ \x^2\lg(t/\x^2)}{2t(\lg t)^2} + O\bigg(\frac{\x^2}{t(\lg t)^2}\bigg).
 \eeq
Finally,  noting  $\int _{t/2 }^{t} \frac{\lg (t-s)}{(\lg s)^2}ds +\int_0^{t/2}\frac{ds}{\lg s} =t/\lg t +O(t/(\lg t)^2)$ and recalling  (\ref{E1}) 
as well as   $F_0(t,\x)-F_0(t,0) = I[D_0(\x^2)] + I[D_>(\x^2)] + I[D_<(\x^2)]$  we conclude the formula of Lemma \ref{lem3.5}. \qed

\v2\n
{\it Proof of Theorem \ref{thm2}}. ~ From Lemmas \ref{lem3.4} and  \ref{lem3.50} it follows that for $1<|\x|<M\sqrt t$,
$$\frac{F(t,\x,1)}{p_t(\x)} = \frac{N(\k t) -\pi p_t(1)+ F_0(t,\x)-F_0(t,0)}{p_t(\x)} +  O(1).$$
By  Lemma  \ref{lem3.5} we can write  the right-hand side as 
$$2\pi tN(\k t)  + 2\pi t\bigg[-N(\k t) +\frac{1}{\lg t}\bigg]\frac{p_t(\x)-p_t(0)}{p_t(\x)} +   \frac{\pi \x^2}{(\lg t)^2}\bigg[\lg\bigg(\frac{t}{\x^2}\bigg) +O(1)\bigg] + O(1), $$
in which the second term  reduces to $O(\x^2/(\lg t)^2)$  since $N(\k t) = 1/\lg t + O(1/(\lg t)^2)$. In view of (\ref{Area}) this yields the formula  
of Theorem \ref{thm2}. \qed

\section{On the asymptotic expansion related to  $N(\la)$}

Our  arguments made in this section are based on the following formula due to C. J.
 Bouwkamp \cite{B}: for  $\la>0, \sigma \geq 0$
\beqn\label{Bouw}
 \int_0^\infty \frac{u^{\sigma-1}e^{-\la u}}{(\lg u)^2 +\pi^2}du = \int_0^s\frac{\sin \pi x}{\pi}\Ga(\sigma +x)e^{-x \lg \la}dx + \frac{\th(\la,s) \Ga(\sigma +s)}{\pi^2 \la^s},
\eeqn
where $s$ may be an  arbitrary positive number and  $|\th(\la,s)|\leq 1$ (Eq (11) of \cite{B}; see  also  Dorning et al  \cite{DNT}  Eq (17)  for the case $\sigma=0$).  From this he obtains the asymptotic expansion of $N(\la)$ in powers of $1/\lg \la$ as $\la\to\infty$. We shall extend the  argument of \cite{B}.  Take $\sigma =0$  in (\ref{Bouw}) so that the integral on the left gives $N(\la)$ and  
substitute from  Euler's relation $\Ga(x)\Ga(1-x)= \pi/\sin \pi x$ as well as $\la= \a t$  to see  that as $t\to\infty$ 
\beqn\label{RF1}
N(\a t) = \int_0^s\frac{1}{\Ga(1-x) \a ^x}e^{-x \lg t}dx + O\bigg(\frac1{t^{s}}\bigg).
\eeqn

\begin{Lem}\label{lem_RF1}
~ For any constant $\a>0$, the function  $N(\a t)$ admits the asymptotic expansion 
$$N(\a t) \sim \frac1{(\lg t)} \sum_{n=0}^\infty  \frac{ a_{n}(\a)}{   (\lg t)^{n}}~~~~~~~ ( t  \to  \infty)$$ 
with the constants  $a_n(\a)$ determined by  
\beqn\label{RF3}
\frac1{\Ga(1-x) \a^x} = \sum_{n=0}^\infty \frac{ a_n(\a)}{n!} x^n~~~~ (|x|<1).
\eeqn
\end{Lem}
\n
\pf ~ The integral on the right-hand side of (\ref{RF1})  is a Laplace transform, in the variable $\lg t$, of a function regular at $x=0$ and it is  standard to infer that the integral admits 
an asymptotic expansion in  the powers of  $1/\lg t$ with the coefficients determined by (\ref{RF3}), showing   the lemma.    \qed

\begin{Lem}\label{lem_RF2}
~ For each constant $\a>0$, as $t\to\infty$
$$\int_0^t N(\a s)ds  \sim \frac{t}{(\lg t)} \sum_{n=0}^\infty   \frac{b_{n}(\a)}{   (\lg t)^{n}},$$ 
$$\int_0^t \Big[-\a s N'(\a s)\Big] ds  \sim \frac{t}{(\lg t)^2} \sum_{n=0}^\infty  \frac{ (n+1)b_{{n}}(\a)}{(\lg t)^{n}}$$ 
with the constants  $b_n(\a)$ determined by  
$$\frac1{(1-x)\Ga(1-x) \a^x} = \sum_{n=0}^\infty \frac{ b_n(\a)}{n!} x^n~~~~ (|x|<1).$$
\end{Lem}
\n\pf~       Integrate the both sides of  (\ref{RF1})  to see that
\begin{eqnarray}\label{RF2}
\int_{1/\a}^t N(\a u)du &=& t\int_0^s\frac{1}{(1-x)\Ga(1-x) \a ^x}e^{-x \lg t}dx\\
 && -  \int_0^s\frac{1/\a}{(1-x)\Ga(1-x)}dx + O\bigg(\frac1{t^{s-1}}\bigg),\nonumber
\end{eqnarray}
which is valid for every $s>0$.  Note that $1/(1-z)\Ga(1-z) =1/\Ga(-z)$ is an entire function.  Since $O(1)$ term is negligible  for the asymptotic expansion, we obtain the first formula as in the preceding proof. One can  verify the second one  in a similar way but by employing (\ref{Bouw}) with $\sigma=1$.  Alternatively, 
note that the left-hand side integral of the second formula equals $ \int_0^t N(\a s)ds -tN(\a t)$. Then, we readily derive it from  the first one combined with the preceding lemma.   \qed
\v2

 From  the Weierstrass product formula 
 $$\frac1{\Ga(1+x)} = e^{\ga x}\prod_{n=1}^\infty e^{-x/n}\bigg(1+\frac{x}{n}\bigg)$$
 (cf. \cite{L})  it follows that if $\zeta(z) = \sum_{k=1}^\infty k^{-z}$, 
\beqn\label{RF0}
\lg \frac1{\Ga(1-x)}= - \,\ga x - \sum_{n=2}^\infty \frac1{n}\zeta(n)x^n,
\eeqn
hence
$$ \sum_{n=0}^\infty \frac{ a_n(\a)}{n!} x^n = \exp\bigg[  -(\ga +\lg \a)x  - \sum_{n=2}^\infty \frac1{n}\zeta(n)x^n\bigg].$$
Similarly,  using $(1-x)^{-1} =\exp{\sum_{n=1}^\infty n^{-1} x^n}$ we obtain from (\ref{RF0}) 
$$ \sum_{n=0}^\infty \frac{ b_n(\a)}{n!} x^n = \exp\bigg[- (\ga +\lg \a -1)x  - \sum_{n=2}^\infty \frac1{n}[\zeta(n)-1]x^n\bigg].$$
From these identities  one can  readily computes the first several terms  of $(a_n(\a))_{n=0}^\infty$ and those of $(b_n(\a))_{n=0}^\infty$ as exhibited in (\ref{Rmk2}) and (\ref{Rmk3}), respectively.
\vskip4mm


\end{document}